\documentclass[12pt]{article}
\usepackage{amsmath,amsfonts,amssymb,amsthm,pstricks,pst-plot,pst-node,eucal}
\usepackage[dvips]{graphics}

\newcommand{\Z}{\mathbb{Z}}
\newcommand{\C}{\mathbb{C}}
\newcommand{\R}{\mathbb{R}}

\newcommand{\fS}{\mathfrak{S}}
\newcommand{\fSt}{\mathfrak{S}_{3D}}

\newcommand{\Pb}[1]{\mathbb{P}_{\textup{#1}}} 
\newcommand{\Ke}[1]{\mathbb{K}_{\textup{#1}}} 
\newpsobject{showgrid}{psgrid}{subgriddiv=1,griddots=5,gridlabels=0pt}
\newcommand{\oh}{\tfrac12}
\newcommand{\LV}{\Lambda^{\frac\infty2}V}
\newcommand{\ul}[1]{{e}_{#1}}
 
\newcommand{\vac}{v_\emptyset}
\newcommand{\la}{\lambda}
\newcommand{\Ex}{\mathbb{E}}

\DeclareMathOperator{\tr}{tr} 
\DeclareMathOperator{\Ad}{Ad} 
\DeclareMathOperator{\Ai}{Ai} 
\DeclareMathOperator{\const}{const} 

\newtheorem{Theorem}{Theorem}

\begin{document}

\title{Symmetric functions and random partitions}
\author{Andrei Okounkov\thanks{
 Department of Mathematics, University of California at
Berkeley, Evans Hall \#3840, 
Berkeley, CA 94720-3840. E-mail: okounkov@math.berkeley.edu.
Partial financial support by NSF grant DMS-0096246, 
Sloan foundation, and the Packard foundation is
gratefully acknowledged. 
}
}
\date{}
\maketitle

\begin{abstract}
These are notes from the lectures I gave at the 
NATO ASI ``Symmetric Functions 2001'' at the 
Isaac Newton Institute in Cambridge (June 25 -- July 6, 2001). 
Their goal is an informal introduction to 
asymptotic combinatorics related to partitions. 
\end{abstract}

These notes are based on the three lectures that I gave in 
in Cambridge in July of 2001. Some parts, especially the
computations, are given here in more detail than 
in the actual lectures; other parts, such as connections
to geometry of the moduli spaces of curves, have been omitted.  
Our goal in these lectures was to give an informal 
introduction to the results contained in \cite{BOO,O2,OR} (there are
also many closely related papers such as \cite{BO,CK,CKP,GTW,J0,J1,J3}), as well
as to the general subject of asymptotic combinatorics
related to partitions. No attempt, however, will be made here to 
survey this already quite large and very rapidly growing subject. 
Many further references to the literature will be pointed out
in text. 

I want to take this opportunity to thank many people who made 
this NATO ASI in Cambridge such a wonderful and memorable meeting: the 
organizers, the sponsors, the lecturers, and, of course, the participants. 
There is also a very long list of people who, in addition to 
their vast influence on the field itself, influenced my 
understanding of it through numerous fruitful discussion,
in particular: P.~Deift, P.~Diaconis, S.~Fomin, A.~Its, K.~Johansson,
 R.~Kenyon,
R.~Stanley, C.~Tracy,  A.~Vershik, H.~Widom, and many others, including,
certainly, my collaborators A.~Borodin, G.~Olshanski, and N.~Reshetikhin.

It is very tragic that one of the pioneers of the field
and, also, one of the  organizers of this meeting, Sergei Kerov, 
did not live to participate in it. Most of the material in these
lectures is directly related to his groundbreaking and lasting work 
and I want to dedicate these lecture to his memory.

\section{Portrait of a large random partition}

\subsection{}

A partition is a nonincreasing sequence 
$$
\lambda= (\lambda_1 \ge \lambda_2 \ge \lambda_3 \ge \dots \ge 0)
$$
of nonnegative integers such that $\lambda_i=0$ for all
sufficiently large $i$. Similarly, a plane partition 
is a matrix of nonnegative integers
with nonincreasing rows and columns and finitely many nonzero
entries, for example
\begin{equation}
  \label{pi}
 \pi = 
\left (\begin {array}{cccc} 5&3&2&1\\\noalign{\medskip}4&3&1&1
\\\noalign{\medskip}3&2&1\\\noalign{\medskip}2&1\end {array}
\right )\,. 
 \end{equation}
A standard geometric image associated to a partition is its
diagram. There are several competing traditions of drawing 
the diagrams of partitions; we follow the one illustrated in
Figure \ref{f2} which shows the 
diagram of the partition $(8,5,4,2,2,1)$. 

\begin{figure}[!h]
  \begin{center}
    \scalebox{0.3}{\includegraphics{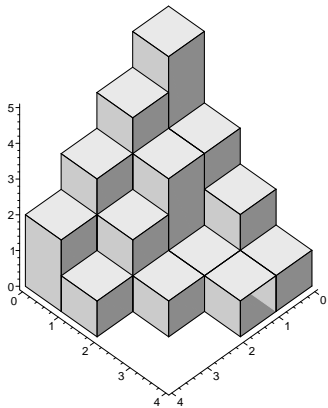}} 
    \caption{A 3-dimensional diagram $\pi$}
    \label{f1}
  \end{center}
\end{figure}

Similarly, the
diagram of a plane partition is a 3-dimensional object 
which looks like a stack of cubes pushed into a corner.
The 3D diagram corresponding to the partition \eqref{pi}
is shown in Figure \ref{f1}

\subsection{} \label{RSK}

Naturally, partitions are inseparable from symmetric
functions. Similarly, application of symmetric
functions to enumeration of plane partitions
is a combinatorial classics, an illustrated account of
which can be found, for example, in \cite{Br}.

Our goal in these lectures is to explain some recent 
progress in applying symmetric function theory to get rather
refined information about 
random partitions, both ordinary and 3-dimensional. 
Random partitions, in addition to being interesting
simply for their own sake, arise in a large variety of other 
contexts. 

One of the most important situations in which random 
partitions arise is the study
of increasing and decreasing subsequences in a random
permutation, see for example \cite{AD,D1}. Namely, the 
celebrated Robinson-Schensted correspondence associates
to any permutation $g\in S(N)$
a partition $\lambda_{RS}(g)$ of $N$, which by Greene--Fomin 
extension (see e.g. \cite{BF}) of Schensted's theorem \cite{Sch} 
encodes the increasing and decreasing
subsequences of $g$. Concretely, the number
$$
\lambda_{RS}(g)_1 + \dots + \lambda_{RS}(g)_k
$$
of squares in the first $k$ rows of $\lambda_{RS}(g)$ equals
the maximal cardinality of a union of $k$ increasing
subsequences of $g$. The same is true about columns
of $\lambda_{RS}$ and decreasing subsequences of $g$. 

Increasing subsequences of random permutations play a central
role, for example, in the analysis of many random growth processes. 
Therefore, one wants to understand  
 the push-forward of the uniform probability
measure on $S(N)$
under the map $g\mapsto \lambda_{RS}(g)$. This push-forward
is well known to be the \emph{Plancherel measure} on 
partitions of $N$ given by
$$
\Pb{Pl}(\lambda)= \frac{(\dim\lambda)^2}{N!} \,,
\quad |\lambda|=N\,, 
$$
where $\dim\lambda$ is the dimension of the corresponding
irreducible representation of the symmetric group. 

The number $\dim\lambda$
has a very transparent combinatorial meaning: it is 
the number of chains of partitions (also known as standard
tableaux) of the form 
\begin{equation}
\lambda = \lambda^{(0)} \searrow \lambda^{(1)} \searrow \lambda^{(2)}
\searrow \dots \searrow \lambda^{(N)}=\emptyset\,, \label{stab}
\end{equation}
where $\mu \searrow \nu$ means that $\nu$ is obtained from 
$\mu$ by removing a square. There are several 
useful formulas for the number $\dim\lambda$, in particular, the 
hook formula. The textbooks \cite{Sa,St} are a wonderful 
place to learn about these things.

\subsection{}

\begin{figure}[!h]
\centering
\scalebox{.7}{\includegraphics{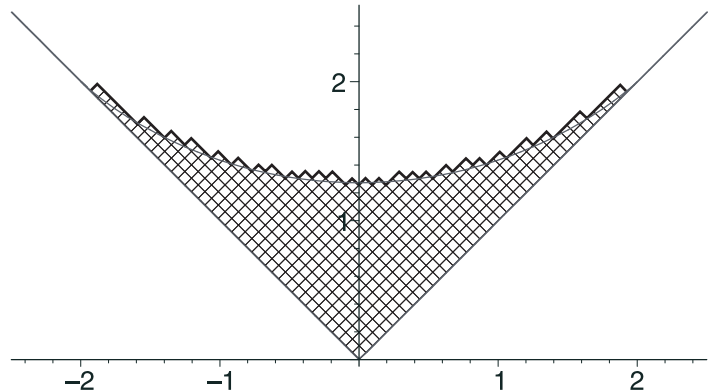}}
\caption{The limit shape of a Plancherel-typical diagram.}
\label{f3}
\end{figure}

The following law of large numbers for the Plancherel measure
was discovered by Logan and Shepp \cite{LS} and, independently,
Vershik and Kerov \cite{VK1,VK2}. Consider a diagram of 
a partition of $N$ as in Figure \ref{f2} and scale it in both
directions by factor of $\sqrt{N}$ (note that the rescaled
diagram has unit area). The rescaled diagram for a 
Plancherel-typical
partition of a large number $N$ will then look 
something like the
diagram in Figure \ref{f3}. 
Now look at the zigzag path in Figure \ref{f3}
formed by the boundary of the diagram. This zigzag is a graph of
a piecewise linear function and, in this way, 
the Plancherel measure becomes a measure on continuous (and even 
Lipschitz with constant $1$) functions.
The discovery of Logan-Shepp-Vershik-Kerov is that,
 as $N\to\infty$, this measure converges  (see \cite{LS,VK1,VK2}
 for precise statements) to the delta measure 
supported on the following function 
\begin{equation*}
\Omega(u)=\begin{cases}
{\frac2\pi\, \left(u \arcsin(u/2) + \sqrt{4-u^2}\right)}\,, & |u|\le 2 \,,\\
|u|\,, & |u|>2 \,.
\end{cases} 
\end{equation*}
This formula looks, at first, rather formidable, but its derivative is much
simpler, namely
\begin{equation}
\Omega'(u) = \frac{2}{\pi} \, \arcsin\frac{u}{2} \,,
\label{Omp}
\end{equation}
naturally extended by $\pm 1$ for $u\notin [-2,2]$. 

In fact, the derivative $\Omega'(u)$ is a very natural object
because it measures how often the zigzag path in Figure \ref{f3}
goes up (respectively, down) in the neighborhood of a given 
point on the boundary 
of the limit shape. In other words, if we think of the zigzag
path as a random sequence of ``ups'' and ``downs'', $\Omega'(u)$  
is related to the probability to observe an ``up'' in a given 
position. Traditionally, probabilities of this kind are
called \emph{correlation functions}. More precisely, the 
the probability to observe ``up'' in a given 
position is the 1-point correlation function of our
random sequence of ``ups'' and ``downs''. 

\subsection{}

It is now natural to ask about higher correlations, for example,
what is the probability to observe several consecutive
``ups'', or, more generally, any given pattern of ``ups'' and
``downs'' ? 

\begin{figure}[!h]
  \begin{center}
    \scalebox{-0.5}[0.5]{\includegraphics{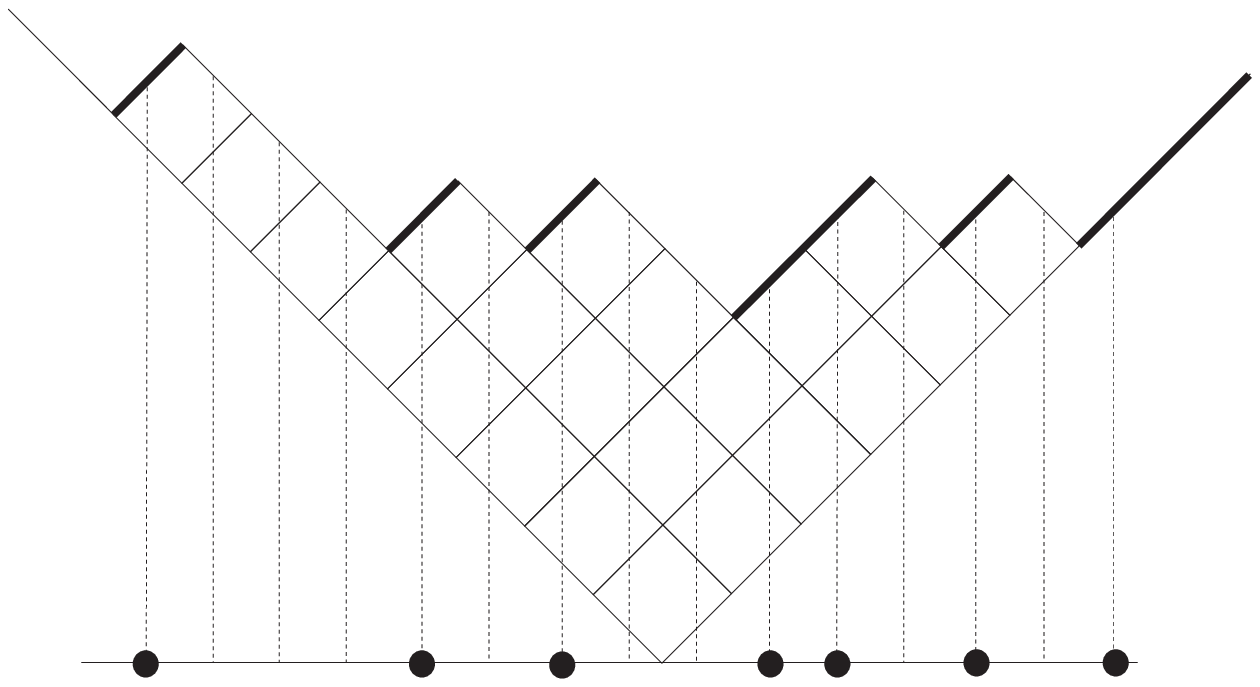}}
  \caption{Geometric meaning of $\fS(\lambda)=\{\bullet\}$}
\label{f2}
\end{center}
\end{figure}

Such questions can be conveniently formalized using the 
following set 
\begin{equation}
\fS(\lambda)=\left\{\lambda_i - i+\tfrac12\right\} \subset \Z + \tfrac12\,.\label{fS}
\end{equation}
This is the set of ``downs'' of the boundary of $\lambda$ as 
illustrated in Figure \ref{f2}.

We remark that the map $\lambda \mapsto \fS(\lambda)$ is a
bijection of the set of partitions and the set of subsets
$S\subset \Z+\frac12$ such that
$$
\left|S\setminus \left(\Z+\tfrac12\right)_{<0}\right| =
\left|\left(\Z+\tfrac12\right)_{<0}\setminus S\right| 
< \infty \,.
$$
In other words, the set $\fS(\lambda)$ has equally
many positive elements and negative ``holes'' (this number
is the length of the diagonal in the diagram $\lambda$). 

We now take a finite set
$$
X=X(N)=\{x_1,\dots,x_n\}\subset \Z+\tfrac12\,,
$$
which is allowed to depend on $N$ in such a way that
$$
\frac{x_i}{\sqrt{N}} \to u_i \in [-2,2]\,, \quad N\to \infty\,,
$$
and all differences $x_i-x_j$ have a limit, finite or 
infinite. The set $X$ is the set where we want to observe
the ``downs''. The parameters $u_i$ describes the global position 
of our observations on the limit shape $\Omega$ while
the differences $x_i-x_j$ record their relative local 
separation. 

It is clear
the following probabilities (a.k.a.\ correlation functions)
$$
\Pb{Pl}\left(X\subset\fS(\lambda)\right) =
\frac1{N!} \sum_{X\subset\fS(\lambda)} (\dim\lambda)^2 
$$
to observe ``downs'' in positions $X=\{x_1,\dots,x_n\}$ include all 
possible information about the local shape of the random 
partition $\lambda$. As it turns out, all these probabilities
have a very pretty $N\to\infty$ limit, namely, we have the 
following:

\begin{Theorem}[\cite{BOO}]\label{T1} For $X$ as above we have
  \begin{equation}
\Pb{Pl}\left(X\subset\fS(\lambda)\right) \to 
\det \Big[\Ke{sin}(x_i-x_j;\phi_i)\Big]_{i,j=1\dots n} \,,
\label{limPl}
\end{equation}
as $N\to\infty$, where 
\begin{equation}
\Ke{sin}(x;\phi)=\frac{\sin \phi x}{\pi x} \,, \quad k \in \Z\,, 
\label{dSin}
\end{equation}
is the discrete sine kernel with parameters
\begin{equation}
\phi_i=\arccos(u_i/2)\,, \quad i=1,\dots,n \label{phiu} \,. 
\end{equation}
\end{Theorem}

In the formula \eqref{dSin}, it is understood that
$$
\Ke{sin}(\infty;\phi) = 0 \,,
$$
which implies that the determinant in \eqref{limPl}
factors into blocks corresponding to distinct
values of $u_i$. In other words, it means that the
occurrences of ``downs'' in different parts of the 
limit shape become independent in the $N\to\infty$ limit,
which, intuitively, is what one expects. 

It is also understood in \eqref{dSin} that
$$
\Ke{sin}(0;\phi) = \frac{\phi}{\pi} \,. 
$$
In particular, the 1-point correlation function, that is,
the density of ``downs'' is
$$
 \Pb{Pl}\left(\{x\}\subset\fS(\lambda)\right) \to
\frac{\arccos u/2}{\pi} \,, \quad \frac{x}{\sqrt{N}}\to u \,. 
$$
This translates into \eqref{Omp} and also explains the role of the 
parameter $\phi$ in the sine kernel: it is the density
parameter. 

Finally, it is natural to extend the function \eqref{phiu}
by $0$ or $\pi$ for $u\notin[-2,2]$. Theorem \ref{T1} remains
valid and says simply that for $u>2$ (resp.\ $u<-2$)
the zigzag in Figure \ref{f3}
goes only up (resp.\ down).

\subsection{}

One encounters a formula which looks just like
the formula \eqref{limPl} in the 
random matrix theory, see e.g.\ \cite{D2,Me}. Concretely,
consider a Gaussian probability measure $\Pb{H}$ on the 
real vector space of Hermitian $N\times N$ matrices
$H$ with density $e^{-\tr H^2}$. The spectrum of 
$\sigma(H)$ of a random matrix $H$ is a random 
$N$-element subset of $\R$. Consider the 
probability density $\rho_{H}(X)$ to have $X\subset \sigma(H)$
for some fixed $X=\{x_1,\dots,x_n\}$, that is,
the probability density to observe eigenvalues
of $H$ in given positions. Formally, this probability
density is defined as follows
$$
\rho_H(x_1,\dots,x_n) \, dx_1 \cdots dx_n = \Pb{H} \left(
\left| \sigma(H) \cap \bigcup_{i=1}^n [x_i,x_i+dx_i]\right|=n\right) \,.
$$
In other words, $\rho_H(x_1,\dots,x_n) \, dx_1 \cdots dx_n$ is 
actually a measure whose values on a given set $A\subset\R^n$ is
related to the probability of finding $n$ eigenvalues of $H$ is $A$.
 
It is important to remember that since $\rho_H(X)$ is a density,
it transforms nontrivially under a change of variables and,
in particular, it scales nontrivially under the scaling of 
variables. 

The $1$-point correlation function $\rho_H(x_1)$
describes the density of eigenvalues of $H$. 
The analog of the limit shape $\Omega$ in this case
is the Wigner semicircle law
$$
\frac1{\sqrt N}\, \rho_H(x) \to \frac{\sqrt{4-u^2}}{2\pi}\,,
\quad \frac{x}{\sqrt{N}}\to u\in[-2,2]\,,\quad N\to \infty \,. 
$$
In other words, scaled by $\sqrt N$, the 
density of eigenvalues converges to the semicircle with 
diameter $[-2,2]$. The normalization $\frac1{\sqrt N}$
is related to the fact that there are $N$ eigenvalues of 
$H$ and, because they all are concentrated
on the segment $[-2\sqrt N, 2\sqrt N]$, the density of 
eigenvalues is expected to be of order $\sqrt N$. 

By the same token, the typical spacing between the eigenvalues
is of order $1/\sqrt{N}$. It is
then meaningful to ask about the  asymptotics of $\rho(X)$ as 
$$
\frac{x_i}{\sqrt{N}} \to u \in [-2,2] \,, \quad
\sqrt{N}(x_i-x_j)\to y_i-y_j \,, 
$$
for some fixed $y_i\in \R$. 
One finds (see \cite{D2,Me}) that in this limit 
\begin{equation}
\frac{\rho(X)}{N^{|X|/2}}\to 
\det \Big[\Ke{sin}(y_i-y_j;r(u))\Big] \,,
\quad r(u)=\sqrt{1-(u/2)^2} \,.\label{limH}
\end{equation}
This looks almost identical to \eqref{limPl}, with  one minor and
one not so minor difference. The minor difference is a different
dependence of the density parameter in the sine kernel on the
global position $u$.
A more important difference is that in \eqref{limPl} the 
variables $x_i-x_j$ take only discrete values, whereas in 
\eqref{limH} we deal with continuous variables. 

In the continuous
case, we can get rid of the density parameter in the sine
kernel by a simple rescaling. After this transformation, the 
global question of eigenvalue density and the local question 
of spacings between the eigenvalues become completely 
decoupled. Much work has been done in the random matrix
theory to show that the eigenvalue spacing have the 
same sine kernel form for many much more general measures
on the space of matrices, in other words, to show that
the occurrence of the sine kernel is \emph{universal}, see
e.g.\ \cite{D2,J2}. 

Comparing, \eqref{limPl} with \eqref{limH}
reveals an striking ``universality of formulas'', namely that 
the formulas are the same even when the meaning is 
different. This parallel, however, is more than 
purely formal. For example, the discrete sine kernel
process shares many features of its continuous 
counterpart, such as the repulsion of eigenvalues
which we will discuss momentarily. Also, 
continuous distributions, such as the distributions
of eigenvalues of a $N \times N$ Gaussian Hermitian
matrix, can be obtained as a suitable limit of 
certain natural measures (of the kind considered
below) on partitions with $N$ parts.

\subsection{} 
It can be inferred from the formula \eqref{limH} that 
the eigenvalues of a random matrix have a strong
tendency to be equally spaced, much more so than 
simply a random collection of points on a line. For
example, the 2-point correlation function (in the
unit density case) is equal to 
$$
\det
\begin{bmatrix}
  1 & \frac{\sin\pi(x-y)}{\pi(x-y)} \\
\frac{\sin\pi(y-x)}{\pi(y-x)} & 1 
\end{bmatrix} = 1 - \left(\frac{\sin\pi(x-y)}{\pi(x-y)}\right)^2 \,,
$$
which is very small (has a double zero) near $x-y=0$ and
has maxima for $x-y=1,2,3,\dots$. 

\begin{figure}[!h]
\centering
\scalebox{.35}[0.25]{\includegraphics{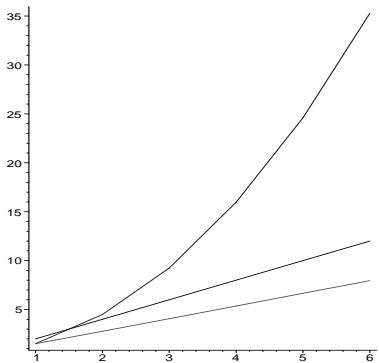}}
\caption{The plot (from top to bottom) of $-\log_2 b_n$, $2n$,
and $-\log_2 a_n$}
\label{fr}
\end{figure}

It was noted by Kerov that a there exists a similar repulsion 
of ``ups'' and ``downs'' in Figure \ref{f3}. Namely, the ``ups'' and
``downs'' tend to be much more uniformly spaced than outcomes of independent
trials, such as, for example, the 
occurrences of ``heads'' in independent coin tosses. 
 Let us illustrate this 
by a numerical example. Suppose that we are in the very 
middle of the limit shape (that is, $u=0$), where the 
density of both ``ups'' and ``downs'' is equal to 
$$
\left.\frac{\arccos(u/2)}{\pi}\right|_{u=0} = \frac12  \,. 
$$
In this case the discrete sine kernel takes an especially
simple form 
$$
\Ke{sin}(x;\tfrac\pi2) = 
\frac{\sin (\pi x/2)}{\pi x} =
\begin{cases}
  \frac12 \,, & x=0\,, \\
  \frac{(-1)^k}{\pi x}\,, & x=2k+1\,, \\
  0\,, & x=\pm 2,\pm 4,\dots \,. \\ 
\end{cases}
$$
Consider the following probabilities
\begin{align*}
  a_n &= \Pb{sin}\left(\textup{up,down,up,down,$\dots$,up,down}\right)\,,\\
  b_n &= \Pb{sin}\left(\textup{up,up,$\dots$,up,down, down,$\dots$,down}\right)
\,, 
\end{align*}
where each string has length $2n$. 
These discrete sine
kernel probabilities are given by certain $(2n)\times(2n)$ determinants. 
If the ``ups'' and ``downs'' were like ``heads'' and ``tails'' of 
independent coin tosses, both $a_n$ and $b_n$ 
would equal $2^{-2n}$. But as the plot of  
$\log_2 a_n$ and $\log_2 b_n$ in Figure \ref{fr} shows
$$
 a_n \gg 2^{-2n} \gg b_n \,,
$$
which means that the discrete sine kernel process very strongly prefers
even spacing of its ``ups'' and ``down''. 

\subsection{}

The connection between the Plancherel measure and random 
matrices becomes exact near the edge of the Logan-Shepp-Vershik-Kerov
limit shape and the Wigner semicircle, respectively. In other
words, the maximal parts $\lambda_1,\lambda_2,\dots$ of 
a Plancherel random partition, suitably scaled, 
 behave in $N\to\infty$ limit
\emph{exactly} like the maximal eigenvalues of an Hermitian
random matrix. For $\lambda_1$ it was proven by Baik, Deift, and
Johansson in \cite{BDJ}, who also conjectured that the same 
holds for any $\lambda_i$, $i=2,3,4,\dots$ (for $i=2$ they 
verified this in \cite{BDJ2}). 

By Schensted's theorem,
mentioned in Section \ref{RSK}, this means that the longest
increasing subsequence of a random permutation is, asymptotically,
distributed like the largest eigenvalue of a Hermitian random 
matrix. This distribution is known as the Tracy-Widom distribution
\cite{TW1} and is given explicitly in term of a certain solution 
of the Painlev\'e II equation. 

The random point process describing the maximal eigenvalues
of a random Hermitian matrix $H$ is known as the \emph{Airy ensemble}
because the correlation function of this process are described
by the same determinantal formula as in \eqref{limH} with 
the sine kernel replaced by the Airy kernel
\begin{equation}
  \label{KAir}
  \Ke{Airy}(x,y)= \frac{\Ai(x) \, \Ai'(y) - \Ai'(x) \, \Ai(y)}{x-y}\,.
\end{equation}
Here $\Ai(x)$ is the Airy function
\begin{equation}
\Ai(x) = \frac1{2\pi} \int_{-\infty}^{\infty} 
e^{is^3/3+ixs} \, ds = \frac1{\pi} \int_{0}^{\infty} 
\cos\left(\frac{s^3}{3} + sx\right)\, ds \,.\label{IA}
\end{equation}
The joint distributions of the largest, 2nd largest, and so on 
points of the Airy ensemble satisfy quite complicated nonlinear
PDEs, generalizing the Painlev\'e II equation for the Tracy-Widom 
distribution. The description of the Airy ensemble 
by its correlation functions is much 
simpler and this is the description that we will use in what
follows.

\subsection{}

The first proof of the full Baik-Deift-Johansson conjecture was
given in \cite{O1} by establishing a direct asymptotics correspondence
between permutations and matrices via a certain correspondence
between counting maps on surfaces
and counting branched covering of the sphere. A different proof was soon
after found, independently, in \cite{BOO} and \cite{J1}. It is this 
second proof that will be explained in these lectures. 

It is clear from Theorem \ref{T1} that
all correlation become trivial at $u=2$, which is because near
$u=2$ the ``downs'' become so rare that they occur infinitely 
far away from each other. In other words, for a Plancherel 
typical partition $\la_i - \la_j \to \infty$ as $N\to\infty$ for
any fixed $i\ne j$. In fact, $\la_i-\la_j$ should be typically
of order $N^{1/6}$, which can be seen heuristically as follows. 

The expected number of $\lambda_i$'s that are larger
than $2\sqrt N - C$ is equal to 
\begin{align*}
  \Ex \left| \fS(\lambda) \cap [2\sqrt N - C, \infty) \right| &\sim
  {\sqrt N}\int^{2}_{2-C/\sqrt N} \frac{\arccos(u/2)}{\pi} \, du \\
  &\sim \frac{2 C^{3/2}}{3 \pi N^{1/4}}
\end{align*}
for $C \ll \sqrt N$. Since there are exactly $k$ parts of $\lambda$ 
that are greater or equal to $\lambda_k$ we expect that for 
$k \ll \sqrt N$
$$
k \approx \const \frac{(2\sqrt N - \lambda_k)^{3/2}}{N^{1/4}} \,,
$$
whence 
$$
2\sqrt N - \lambda_k \approx \const  k^{2/3} \, N^{1/6}
$$
which precisely means that 
$\lambda_i-\lambda_j$ should be typically of order $N^{1/6}$
for large $N$. 

This heuristics explains the scaling in 
the following result (conjectured by Baik, Deift, and Johansson), which is  
the analog of Theorem \ref{T1} for the edge $u=2$ of the
limit shape $\Omega$. 

\begin{Theorem}[\cite{O1,BOO,J1}]\label{TA} As $N\to\infty$, the 
random variables
$$
\frac{\lambda_i-2\sqrt N}{N^{1/6}}\,, \quad i=1,2,3,\dots\,,
$$
converge in joint distribution to the Airy ensemble. 
\end{Theorem}

\subsection{}

We now move on to plane partitions and 3D diagrams. On this 
set, we consider the following measure:
$$
\Pb{3D}(\pi) = \frac1{Z_{3D}} \, q^{|\pi|} \,, \quad q\in [0,1)\,, 
$$
where $|\pi|$ is the volume of a 3D diagram $\pi$ and
$Z_{3D}$ is the normalizing factor
$$
Z_{3D} = \sum_{\pi} q^{|\pi|} = \prod_{n>0} (1-q^n)^{-n}\,,
$$
the second equality being the classical identity due to McMahon. 

In particular, this measure gives equal weight $q^N$ to all 
partitions of size $N$. Is is natural to ask why we consider
this essentially uniform measure and not some analog of the 
Plancherel measure on 3D diagrams (which could involve, 
for example, the number of ways to construct $\pi$ starting
from the empty partition by adding a square at a time). 
The answer is that this measure \emph{is} the proper
analog of the Plancherel measure for the following reason. 
By definition and \eqref{stab}, $\Pb{Pl}(\lambda)$ is proportional
to the number of chains of the form
\begin{equation}
\emptyset=\lambda^{(-N)} \nearrow \dots \nearrow \lambda^{(-1)}
\nearrow \lambda^{(0)} \searrow \lambda^{(1)} 
\searrow \dots \searrow \lambda^{(N)}=\emptyset \,. \label{sstab}
\end{equation}
such that
$$
\lambda^{(0)}= \lambda \,.
$$
If one assembles all these $\lambda^{(i)}$ into an array of partitions,
one sees a structure closely resembling a 3D diagram and, by
construction, the Plancherel measure is the distribution of the 
central slice of this structure induced by the uniform measure
on all chains of the form \eqref{sstab}. 

\begin{figure}[!h]\psset{unit=0.5cm}
  \begin{center}
    \begin{pspicture}(-7,-4)(7,5)
\scriptsize
\showgrid
\psset{dimen=middle}
\psaxes[axesstyle=frame,Ox=-7,Oy=-4,Dx=2,Dy=2,ticks=none](-7,-4)(7,5)
\psdiamond[fillstyle=solid,fillcolor=lightgray](0,4.5)(1,.5)
\psdiamond[fillstyle=solid,fillcolor=lightgray](1,2.0)(1,.5)
\psdiamond[fillstyle=solid,fillcolor=lightgray](2, .5)(1,.5)
\psdiamond[fillstyle=solid,fillcolor=lightgray](3,-1.0)(1,.5)
\psdiamond[fillstyle=solid,fillcolor=lightgray](4,-2.5)(1,.5)
\psdiamond[fillstyle=solid,fillcolor=lightgray](5,-3.0)(1,.5)
\psdiamond[fillstyle=solid,fillcolor=lightgray](6,-3.5)(1,.5)
\psdiamond[fillstyle=solid,fillcolor=lightgray](-1,3.0)(1,.5)
\psdiamond[fillstyle=solid,fillcolor=lightgray](0,1.5)(1,.5)
\psdiamond[fillstyle=solid,fillcolor=lightgray](1,-1.0)(1,.5)
\psdiamond[fillstyle=solid,fillcolor=lightgray](2,-1.5)(1,.5)
\psdiamond[fillstyle=solid,fillcolor=lightgray](3,-3.0)(1,.5)
\psdiamond[fillstyle=solid,fillcolor=lightgray](4,-3.5)(1,.5)
\psdiamond[fillstyle=solid,fillcolor=lightgray](-2,1.5)(1,.5)
\psdiamond[fillstyle=solid,fillcolor=lightgray](-1,0.0)(1,.5)
\psdiamond[fillstyle=solid,fillcolor=lightgray](0,-1.5)(1,.5)
\psdiamond[fillstyle=solid,fillcolor=lightgray](1,-3.0)(1,.5)
\psdiamond[fillstyle=solid,fillcolor=lightgray](2,-3.5)(1,.5)
\psdiamond[fillstyle=solid,fillcolor=lightgray](-3,0.0)(1,.5)
\psdiamond[fillstyle=solid,fillcolor=lightgray](-2,-1.5)(1,.5)
\psdiamond[fillstyle=solid,fillcolor=lightgray](-1,-3.0)(1,.5)
\psdiamond[fillstyle=solid,fillcolor=lightgray](0,-3.5)(1,.5)
\psdiamond[fillstyle=solid,fillcolor=lightgray](-4,-2.5)(1,.5)
\psdiamond[fillstyle=solid,fillcolor=lightgray](-3,-3.0)(1,.5)
\psdiamond[fillstyle=solid,fillcolor=lightgray](-2,-3.5)(1,.5)
\psdiamond[fillstyle=solid,fillcolor=lightgray](-5,-3.0)(1,.5)
\psdiamond[fillstyle=solid,fillcolor=lightgray](-4,-3.5)(1,.5)
\psdiamond[fillstyle=solid,fillcolor=lightgray](-6,-3.5)(1,.5)
\end{pspicture}
    \caption{Projections of horizontal faces of the diagram from 
Figure \ref{f1}}
    \label{f4}
  \end{center}
\end{figure}
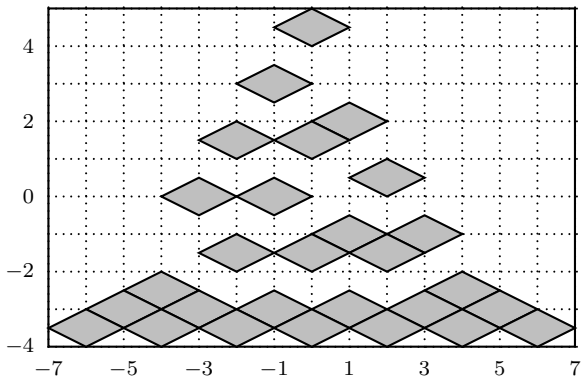

This also suggests that the proper way to look at 3D diagram is
to separate them into their diagonal slices (and also suggest a 
common generalization of the measures $\Pb{Pl}$ and $\Pb{3D}$ to 
be introduced later). In particular, the diagonal slices of the
partition from \eqref{pi} are 
$$
(2) \prec (3,1) \prec (4,2) \prec (5,3,1) \succ (3,1)\succ (2,1)\succ (1)\,,
$$
where $\lambda \succ \mu$ means that $\lambda$ and $\mu$ interlace
(this relation replaces the relation $\lambda\searrow\mu$ in \eqref{sstab}). 
The proper modification of the 
definition \eqref{fS} for the 3D case is 
\begin{equation}
\fSt(\pi)=\left\{\left(j-i,\pi_{ij} -\tfrac{i+j-1}{2}\right)\right
\} \subset \Z\oplus \tfrac12 \Z\,, \label{fSt}
\end{equation}
which has the following geometric meaning. 
The set $\fSt(\pi)\subset \R^2$ is formed by 
the centers of the horizontal faces of a diagram $\pi$ under 
the following mapping 
$$
\R^3 \owns (x,y,z) \mapsto (y-x,z-\tfrac{x+y}{2})\in \R^2\,.
$$
This is illustrated in Figure \ref{f4} which shows the projections
of the horizontal faces of the 3D diagram from Figure \ref{f1}. 

\subsection{}
The analog of the $N\to\infty$ limit for the Plancherel measure
is the $q\to 1$ limit of $\Pb{3D}$. Indeed, as $q\to 1$, larger
and larger diagrams become important. It can be shown (see e.g.\ \cite{OR})
that as $q\to 1$
$$
r^3 |\pi| \to 2\zeta(3)
$$
in probability, where 
$$
r = - \ln q \,.
$$
is a convenient parameter that goes to $+0$ as $q\to 1$. 

Thus, it is very natural to ask what is the asymptotics
of the correlation functions
$$
\Pb{3D}\left(X\subset\fSt(\pi)\right)\,, \quad
X={(t_i,h_i)}\subset \Z\oplus\tfrac12\Z\,,
$$
as $q\to 1$ and the set $X$ is allowed to vary in such a 
way that
$$
r (t_i,h_i) \to (u,v) \in \R^2\,,
$$
while all differences 
$$
(t_i,h_i)-(t_j,h_j) \in \Z\oplus\tfrac12\Z
$$
have a limit, finite or infinite. 

At this point, it should
not come as a big surprise that there is an analog 
of Theorem \ref{T1} for 3D diagrams that answers this
question. To state it, we need the following  $2$-dimensional analog of
the sine kernel. One way to write the sine kernel is like
this:
$$
\Ke{sin}(k;\phi)=\frac{\sin k\phi}{\pi k} = \frac{1}{2\pi i}
\int_{e^{-i\phi}}^{e^{i\phi}} \frac{dw}{w^{k+1}} \,.
$$
A $2$-dimensional analog of this is the following 
\emph{incomplete beta kernel}
\begin{equation}
\Ke{beta}(t,h;z) = \frac1{2\pi i} \int_{\bar z}^z (1-w)^t \, 
\frac{dw}{w^{h+t/2+1}} \,,\label{Kbeta}
\end{equation}
where $z\in\C$ is a parameter that replaces the density 
parameter $\phi\in\R$. The path of the integration in
\eqref{Kbeta} crosses the segment $(0,1)$ for $t\ge 0$ and the 
segment $(-\infty,0)$ for $t<0$. 

Now the analog of Theorem \ref{T1} for 3D diagrams is
the following. By symmetry, we can assume, without
loss of generality that $u\ge 0$. 

\begin{Theorem}[\cite{OR}]\label{T2} For $X$ as above we have
\begin{equation}
\Pb{3D}\left(X\subset\fSt(\pi)\right) \to 
\det \Big[\Ke{beta}(t_i-t_j,h_i-h_j;z)\Big]_{i,j=1\dots n} \,,
\label{lim3D}
\end{equation}
as $q\to 1$, where the parameter $z$ is the point of 
intersection of two circles
$$
\{z, \bar z\} = \{|z|=e^{-u/2}\} \cap 
\{|z-1|=e^{-u/4-v/2}\}\,, \quad \Im z >0\,,
$$
or else, in case these circles do not intersect,
$z=\bar z$ is  the point of $\{|z|=e^{-u/2}\}$ which is the closest to 
$\{|z-1|=e^{-u/4-v/2}\}$. 
\end{Theorem}

\subsection{}

In particular, the limit of the 1-point 
correlation function, that is, the limiting
density of the horizontal faces, is given by 
\begin{equation}
\Ke{beta}(0,0;z)=\frac{\arg z}{\pi} 
= \frac1\pi{\arccos\left(\cosh\frac u2 - \frac{e^{-v}}{2}\right)}\,.
\label{3Dd}
\end{equation}
Here the second equality is literally true if the 
two circles in Theorem \ref{T2} intersect. Otherwise, 
the arccosine has to be extended by $0$ or $\pi$ which,
just like for the Plancherel measure, 
corresponds to trivial correlations. In other words,
when the two circles in Theorem \ref{T2} do not
intersect it means that we are trying to 
measure correlations outside of the limit
shape. 

\begin{figure}[!h]
  \begin{center}
   \scalebox{0.6}{\includegraphics{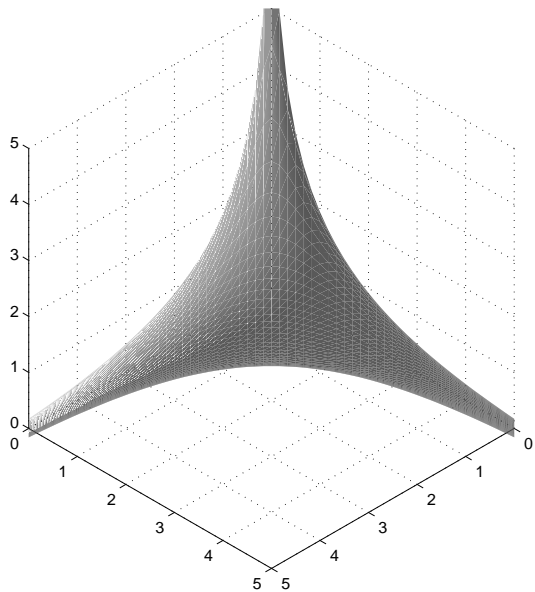}}     
    \caption{The limit shape for 3D diagrams}\label{f5}
     \end{center}
\end{figure}

This limit shape, the existence of which was first
established by Vershik \cite{V} and which was 
described explicitly by Cerf and Kenyon in \cite{CK},
can be reconstructed by integrating the density \eqref{3Dd}.
This integral can be evaluated in terms of the dilogarithm
function. The plot of the limit shape can be seen in 
Figure \ref{f5}.

\section{Schur process and its correlations}

\subsection{}

It is difficult to deny that the best way to solve an 
asymptotic problem is to have a nice exact formula valid 
before the limit, which would make passing to the limit
conceptually clear, if not completely straightforward. A wonderful
feature of our problems is that such exact formulas
are indeed available. 

In fact, exact formulas for correlation functions are
available for much more general class of measures. These
measures, which we call \emph{Schur process}, 
will be described presently. By definition \cite{OR}, the 
Schur process is a measure on sequences of partitions
$\{\lambda(t)\}$, $t\in\Z$, such that $\lambda(t)=\emptyset$
for all but finitely many $t$. The probability 
of a configuration $\{\lambda(t)\}$ is
\begin{equation}
\Pb{Schur}(\{\lambda(t)\}) = \frac1{Z} 
\prod_{m\in\Z+\frac12} 
s_{\lambda(m\pm\frac12)/\lambda(m\mp\frac12)}(A_m)\,,
\label{PSchur}
\end{equation}
where $s_{\lambda/\mu}$ is a skew Schur function and 
$A_m$ is some set of variables. The factor $Z$ in \eqref{PSchur}
in introduced, as usual, for normalization, that is, to make
$\Pb{Schur}$ a probability measure. 
The alphabets $A_m$ and
the choice of signs in \eqref{PSchur} are the 
parameters of the process. In our examples, the 
choice of signs will always be understood to be such that
$$
  \dots  \subset \lambda(-2) \subset  \lambda(-1)
 \subset      \lambda(0) \supset \lambda(1)  \supset \lambda(2) 
\supset \dots \,.
$$
In other words, we would like our partitions $\lambda(t)$ to grow
up to the middle and then decay, just like in \eqref{sstab}.
The general formalism of Schur processes, however, works equally
well for any choice of signs. The above definition of a
Schur process is slightly different from the one given 
in \cite{OR}, which avoids the awkward necessity of choosing
signs, but is essentially equivalent. 

\subsection{}

For example, let
$$
A_m = \{q^{|m|}\}, \quad m\in\Z+\tfrac12\,,
$$
that is, let the alphabet $A_m$ consist of a single 
variable $q^{|m|}$, where $q$ is a parameter 
(e.g.\  a complex number such that $|q|<1$). 
It is a basic fact that
$$
s_{\lambda/\mu}(t,0,0,\dots)=
\begin{cases}
  t^{|\lambda/\mu|}\,, & \lambda \succ \mu\,,\\
0 \,, & \textup{otherwise}\,,
\end{cases}
$$
where $\lambda \succ \mu$ means that $\lambda$ and $\mu$ interlace,
which can be also expressed as $\lambda/\mu$ being a horizontal strip. 
Therefore, in this case the Schur process is supported on the
sequences of the form 
$$
  \dots  \prec \lambda(-2) \prec  \lambda(-1)
 \prec      \lambda(0) \succ \lambda(1)  \succ \lambda(2) 
\succ \dots 
$$
and the probability of such a configuration is 
proportional to 
$$
q^{\sum_{m\in\Z+\frac12} |m| |\lambda(m\pm\frac12)/\lambda(m\mp\frac12)|}
= q^{\sum |\lambda(m)|} \,.
$$
We recognize the measure $\Pb{3D}$, where the diagram 
$\pi$ is given by the sequence $\{\lambda(m)\}$ of its
diagonal slices. 

The general formalism of \cite{O2,OR}, which will be explained
on an example later, gives the following exact formula 
for the correlation functions of the measure $\Pb{3D}$.
We have
\begin{equation}
  \label{3De}
  \Pb{3D}\left(X\subset\fSt(\pi)\right) =
\det \Big[\Ke{3D}((t_i,h_i),(t_j,h_j);q)\Big]_{i,j=1\dots n}\,,
\end{equation}
where the kernel $\Ke{3D}$ is given by the following
contour integral 
\begin{multline}
  \label{K3D}
  \Ke{3D}((t_i,h_i),(t_j,h_j);q) = \\
\frac1{(2\pi i)^2} \int_{|z|=1\pm\epsilon}
\int_{|w|=1\mp\epsilon} \frac1{z-w} \, 
\frac{\Phi_{3D}(t_1,z)}{\Phi_{3D}(t_2,w)}\,
\frac{dz\, dw}{z^{h_1+\frac{|t_1|}2+\frac12}
 w^{-h_2 - \frac{|t_2|}2+\frac12}}\,. 
\end{multline}
In this integral, one choses the top sign if $t_1\ge t_2$ and 
the bottom sign otherwise. The 
function $\Phi_{3D}$ here is the following product 
$$
\Phi_{3D}(t,z) = \frac
{\prod_{m>\max(0,-t)} (1-q^m/z)}
{\prod_{m>\max(0,t)} (1-q^m z)}\,, \quad m\in\Z+\tfrac12 \,.
$$
At first, the integral in \eqref{K3D} looks a little scary,
but in fact it is not and, in particular, its asymptotic
analysis, leading to Theorem \ref{T2} is more or less 
straightforward saddle point analysis, see \cite{OR} and below.

\subsection{}

As another example, take the following specialization of
the Schur process. Let $A_m$ be empty for $m\ne\pm\frac12$
and for $m=\pm\frac12$ take the following specialization of the 
algebra of the symmetric functions
$$
p_1 \mapsto \sqrt\xi\,, \quad p_k\mapsto 0\,, \quad k>1\,,
$$
where $\xi>0$ is a parameter. 
It is well known that this has the following effect
$$
s_\lambda\left(A_{\pm\frac12}\right) = \xi^{|\lambda|/2} \, 
\frac{\dim\lambda}{|\lambda|!} \,.
$$
Thus, in this specialization, the Schur process becomes
a measure on sequences
$$
 \dots  \subset \emptyset \subset \emptyset 
 \subset      \lambda \supset \emptyset \supset \emptyset
\supset \dots 
$$
where the distribution of the unique nontrivial 
partition in the middle is given by
\begin{equation}
  \label{PP}
  \Pb{PP}(\lambda)= \frac1{Z_{PP}} \, \xi^{|\lambda|} \, 
\left(\frac{\dim\lambda}{|\lambda|!}\right)^2\,, \quad 
Z_{PP}=e^\xi\,.
\end{equation}
This measure is the so-called Poissonization of the
Plancherel measure, that is, the mixture of Plancherel
measures on partitions of $N$ for different values of $N$ by means of
the Poisson distribution on $N$ with parameter $\xi$
\begin{equation}
\Pb{Poisson}(N) = \frac1{e^{\xi}} \, \frac{\xi^N}{N!} \,,
\quad N=0,1,2,\dots \,.
\label{Poi}
\end{equation}
In particular, from \eqref{Poi} we can infer the value of the normalization 
constant in \eqref{PP}. 

It is a well known principle (known under various names
such as e.g.\ \emph{equivalence of ensembles}) that the 
$N\to\infty$ asymptotics of the Plancherel measures
should be equivalent to the $\xi\to\infty$ limit of 
their Poissonizations \eqref{PP}. Essentially, this
is because the distribution \eqref{Poi} has mean $\xi$ and 
standard deviation $\sqrt\xi$, so for $\xi\sim N\to\infty$
it is concentrated around $N$.

\subsection{}\label{sKB}
The general formula of \cite{O2,OR} specializes to the
following exact formula for the correlations of the
measure $\Pb{PP}$ (which is a limit case of the main formula of
\cite{BO} and can be found as stated in \cite{BOO} and \cite{J1}) 
\begin{equation}
\Pb{PP}\left(X\subset\fS(\lambda)\right) = 
\det \Big[\Ke{Bessel}(x_i,x_j;\xi)\Big]_{i,j=1\dots n} \,,
\label{PPc}
\end{equation}
where $\Ke{Bessel}$ is the discrete Bessel kernel 
given by
\begin{align}
  \label{KB}
  \Ke{Bessel}(x,y;\xi) & = 
\frac1{(2\pi i)^2} 
\iint_{|z|>|w|} \frac{e^{\sqrt{\xi}(z-z^{-1}-w+w^{-1})}}{z-w} \,
\frac{dz\, dw}{z^{x+\frac12}\, w^{-y+\frac12}} \\
&= \sqrt\xi \, \frac{J_{x-\frac12} \, J_{y+\frac12} -
J_{x+\frac12} \, J_{y-\frac12}}{x-y} \,. \label{KB2}
\end{align}
Here $J_n = J_n(2\sqrt\xi)$ is the Bessel function of
integral order $n$ (recall that $x,y\in\Z+\frac12$) defined
as the coefficient in the following Laurent series
$$
e^{\frac x2(z-z^{-1})} = \sum_{n\in\Z} J_n(x) \, z^n \,.
$$
The equivalence of \eqref{KB} and \eqref{KB2} can be seen
as follows. We have
$$
\frac{1}{z^{x+\frac12}\, w^{-y+\frac12}}
   = - \frac1{x-y}\,\left(z\frac{\partial}{\partial z} + 
w\frac{\partial}{\partial w} + 1 \right) \, 
\frac{1}{z^{x+\frac12}\, w^{-y+\frac12}} \,,
$$
which we can substitute in \eqref{KB} and then integrate
by parts using
$$
\left(z\frac{\partial}{\partial z} + 
w\frac{\partial}{\partial w} + 1 \right) \,
\frac{e^{\sqrt{\xi}(z-z^{-1}-w+w^{-1})}}{z-w} =
\sqrt{\xi} \, \left(1-\frac1{zw}\right)\,
e^{\sqrt{\xi}(z-z^{-1}-w+w^{-1})} \,.
$$
This yields \eqref{KB2}\,.

\subsection{}

We will now explain, on the example of the formula \eqref{PPc},
the general formalism for obtaining such formulas developed 
in \cite{O2,OR}. Consider the following operator in the 
algebra of symmetric functions
$$
\Gamma \cdot f = e^{\sqrt\xi \, p_1} f \,.
$$
Since
\begin{equation}
p_1 \, s_\mu = \sum_{\lambda\searrow\mu} s_\lambda\,,
\label{Pie}
\end{equation}
we conclude that
$$
\Gamma \cdot 1 = \sum_\lambda \xi^{|\lambda|/2} \, 
\frac{\dim\lambda}{|\lambda|!} \, s_\lambda \,.
$$
Consider the following projection operator
\begin{equation}
\delta_X \cdot s_\lambda =
\begin{cases}
  s_\lambda\,, & X\subset\fS(\lambda)\\
0\,, & \textup{otherwise}\,.
\end{cases}\label{dX}
\end{equation}
Then, by construction
\begin{equation}
\Pb{PP}\left(X\subset\fS(\lambda)\right) =
\frac{(\Gamma^* \, \delta_X \, \Gamma \cdot 1, 1)}
{(\Gamma^* \, \Gamma \cdot 1, 1)}\,,
\label{PPc1}
\end{equation}
where we use the standard Hall inner product
on the algebra of symmetric functions
$$
(s_\lambda,s_\mu)=\delta_{\lambda,\mu} \,,
$$
and where 
$$
\Gamma^* = \exp\left(\sqrt{\xi}\,\frac{\partial}{\partial p_1}\right) \,,
$$
is the adjoint of $\Gamma$ with respect
to this inner product. The relation
$$
\left[  \frac{\partial}{\partial p_1}, p_1\right] =1
$$
implies that
\begin{equation}
\Gamma^* \, \Gamma =  e^\xi \, \Gamma \, \Gamma^*\,,
\label{Gc}
\end{equation}
which together with
\begin{equation}
\Gamma^* \cdot 1 = 1
\label{Gf}
\end{equation}
gives another computation of the normalization
constant
$$
Z_{PP} = (\Gamma^* \, \Gamma \cdot 1, 1) = e^\xi \,.
$$
Our goal is to apply a similar trick to evaluate
the numerator in \eqref{PPc1}. For this we need
a good realization of the operator \eqref{dX}. 
This is naturally provided by the fermionic Fock space 
formalism (see, e.g.\ \cite{K}) which will be discussed presently.

\subsection{}
Let the vector space $V$ be spanned by $\ul{k}$, $k\in\Z+\oh$.
The infinite wedge space $\LV$  is, by definition, spanned by vectors 
$$
v_S=\ul{s_1} \wedge \ul{s_2} \wedge  \ul{s_3} \wedge  \dots\,,
$$
where $S=\{s_1>s_2>\dots\}\subset \Z+\oh$ is such a subset that
both sets
$$
S_+ = S \setminus \left(\Z_{\le 0} - \oh\right) \,, \quad
S_- = \left(\Z_{\le 0} - \oh\right) \setminus S 
$$
are finite. We equip $\LV$ with the inner product 
in which the basis $\{v_S\}$ is orthonormal. In particular,
we have the vectors
$$
v_\lambda=\ul{\lambda_1-\frac12} \wedge \ul{\lambda_2-\frac32}
\wedge \ul{\lambda_3-\frac52} \wedge \dots \,,
$$
where $\lambda$ is a partition. 
The vector 
$$
\vac = \ul{-\frac12} \wedge \ul{-\frac32} \wedge \ul{-\frac52} \wedge \dots
$$
is called the vacuum vector. 

By construction, the $\C$-linear map 
$$
s_\lambda \to v_\lambda
$$
is an isometric embedding of the algebra of symmetric functions
into $\LV$. It follows from \eqref{Pie} that
the operator of multiplication by $p_1$ is given  in 
this realization by the natural action on $\LV$ of the following
matrix 
$$
p_1 \mapsto \alpha = \big[\delta_{k,l-1}\big]_{k,l\in\Z+\frac12} =
\begin{bmatrix}
\ddots & \ddots \\  
&0 & 1 \\
&& 0 & 1 \\
&&& 0 & 1 \\
&&&& 0 & 1 \\
&&&&& 0 & 1 \\
&&&&&& \ddots & \ddots 
\end{bmatrix} \in \mathfrak{gl}(V) \,.
$$
As we shall see in a moment, the operator \eqref{dX} is
given in this realization by something equally natural. 

Consider the action of the matrix unit
$$
E_{ii} = \big[\delta_{k,i}\delta_{l,i}\big]_{k,l\in\Z+\frac12} =
\begin{bmatrix}
\ddots \\  
&0 &  \\
&& 1 &  \\
&&& 0 &  \\
&&&&& \ddots  
\end{bmatrix} \in \mathfrak{gl}(V) \,.
$$
It is clear that
$$
E_{ii} \, v_\lambda = 
\begin{cases}
  v_\lambda\,, & i\in\fS(\lambda)\\
0\,, & \textup{otherwise} \,,
\end{cases}
$$
and hence
$$
\delta_X = \prod_{x\in X} E_{xx}  \,.
$$
This is already very nice, but we can actually do better
by introducing the fermionic operators $\psi_i$ which 
are, in some sense, the square roots of the operators $E_{ii}$.

By definition, the operator $\psi_i$ is the wedging with the
vector $\ul{i}$
$$
\psi_i \, v = \ul{i} \wedge v \,,
$$
and the operator $\psi^*_i$ is the adjoint operator with
respect to the inner product on $\LV$. One checks from
definitions that
$$
E_{ii} = \psi_i \, \psi^*_i \,,
$$
and hence the formula \eqref{PPc1} becomes
\begin{equation}
\Pb{PP}\left(X\subset\fS(\lambda)\right) =
\frac1{e^\xi}
\left(e^{\sqrt\xi \,\alpha^*} \, \prod_{x\in X} \psi_x\, \psi^*_x \,
e^{\sqrt\xi \, \alpha}  \cdot \vac, \vac\right)\,.
\label{PPc2}
\end{equation}

\subsection{}

Now we want to exploit \eqref{Gc} and \eqref{Gf}. 
By definition, set
\begin{equation}
\Psi_x = \Ad\left(e^{\sqrt\xi \,\alpha^*} \, e^{- \sqrt\xi \, \alpha}
\right) \cdot \psi_x \,,
\label{Psix}
\end{equation}
and similarly define $\Psi_x^*$, 
where $\Ad$ denotes the action by conjugation. Note that since the 
operators $\alpha$ and $\alpha^*$ commute up to a central 
element, the ordering of these operators inside the $\Ad$ sign 
is immaterial.  Combining \eqref{Gc}, \eqref{Gf},  and \eqref{Psix}
with \eqref{PPc2} gives the following
\begin{equation}
\Pb{PP}\left(X\subset\fS(\lambda)\right) =
\left( \prod_{x\in X} \Psi_x\, \Psi^*_x \,
  \cdot \vac, \vac\right)\,.
\label{PPc3}
\end{equation}
It is now a pleasant combinatorial exercise, known as the 
Wick lemma, to verify that
$$
\left( \prod_{x\in X} \Psi_x\, \Psi^*_x \,
  \cdot \vac, \vac\right) = 
\det \left[ \left(\Psi_x\, \Psi^*_y 
  \cdot \vac, \vac\right) \right]_{x,y\in X}\,.
$$
This reduces the verification of \eqref{PPc} to proving
that
\begin{equation}
\left(\Psi_x\, \Psi^*_y 
  \cdot \vac, \vac\right) = \Ke{Bessel}(x,y;\xi) \,.
\label{Kk}
\end{equation}
Since \eqref{KB} gives, essentially, the 
generation function for the kernel $\Ke{Bessel}$, it
is natural to introduce similar generation functions
\begin{align*}
  \psi(z) & = \sum_{k\in\Z+\frac12} z^k \, \psi_k\,, \\
\psi^*(w) & = \sum_{k\in\Z+\frac12} w^{-k} \, \psi^*_k\,. 
\end{align*}
One then checks directly that
$$
\Psi(z) = \Ad\left(e^{\sqrt\xi \,\alpha^*} \, e^{- \sqrt\xi \, \alpha}
\right) \cdot \psi(z)  = 
e^{\sqrt{\xi}(z-z^{-1})} \, \psi(z) \,,
$$
which together with
$$
\left(\psi(z)\, \psi^*(w)
  \cdot \vac, \vac\right) = \sum_{k=-\frac12,-\frac32,-\frac52,\dots}
z^k w^{-k} = \frac{\sqrt{zw}}{z-w}\,, \quad |z|>|w|
$$
implies that
$$
\left(\Psi(z)\, \Psi^*(w)
  \cdot \vac, \vac\right) = 
\frac{\sqrt{zw}}{z-w} \, e^{\sqrt{\xi}(z-z^{-1}-w+w^{-1})}\,,
$$
which proves \eqref{Kk} and completes the proof of \eqref{PPc}.

\section{Asymptotics}

\subsection{}

The integral representations of the type \eqref{K3D}
and \eqref{KB} are particularly suitable for 
asymptotic analysis by the classical method 
of steepest descent (see e.g.\ \cite{BeO,Ol}). The essence of this method
is that the asymptotics of the integral
\begin{equation}
\int e^{M S(x)} \, dx\,, \quad M\to\infty\,,
\label{IM}
\end{equation}
for a real function $S(x)$ is dominated by
arbitrarily small neighborhoods of 
the points where the function $S(x)$ takes
its maximal value. If $x_0$ is such a point
then $S'(x_0)=0$ and hence 
$$
S(x)= S(x_0) + \frac{S''(x_0)}{2} (x-x_0)^2 + O((x-x_0)^3)\,.
$$
Assume that $S''(x_0)\ne 0$ (and hence $S''(x_0)< 0$) and introduce
a new $t$ variable by 
$$
{x-x_0} = \frac{t}{\sqrt{- M S''(x_0)}} \,.
$$
Then for fixed $t$ 
$$
M S(x) = M S(x_0) - \frac{t^2}{2} + o(1) 
$$
and we see that the leading contribution of $x_0$ to \eqref{IM}
is given by
$$
e^{M S(x_0)}
 \int_{-\infty}^{\infty} e^{-t^2/2} dx = 
\sqrt{2\pi} \, \frac{e^{M S(x_0)}}{\sqrt{-M S''(x_0)}} \,.
$$
Similar considerations apply in the case when the maximum
of $S(x)$ at $x_0$ is degenerate, that is, when $S''(x_0)=0$,
but, for example, $S^{IV}(x_0)<0$. 

For similar integrals in the complex plane
$$
\int_\gamma e^{M S(z)} \, dz\,, \quad M\to\infty\,,
$$
where $\gamma$ is some contour, one first deforms the
contour $\gamma$ to make is pass through the critical 
points $z_0$  
of $S(z)$ (also known as saddle points) in such a way that the real 
part of $S(z)$ has a local maximum along $\gamma$ at $z_0$.
After that, one proceeds as with the integral \eqref{IM}. 

A qualitative conclusion from this argument is that 
\emph{such integrals never have finite limit as $M\to\infty$}.
In other words, whenever one expects some finite limit,
such as, for example, we expect \eqref{dSin}
to be the $\xi\to\infty$ asymptotics of \eqref{KB},
the contributions of saddle points computed by the 
method above will be simply zero. This makes one wonder
how one can ever get some finite limit, because, after all,
\eqref{dSin} is the $\xi\to\infty$ asymptotics of \eqref{KB}
as we will see in a moment. The answer is that, as we 
deform the contours as required by the steepest descent
method, we will be picking up residues, and those will 
have a nontrivial finite limit.

\subsection{}
\begin{figure}[!h]
\begin{center}
\begin{pspicture}(-2.2000000,-2.2000000)(2.200000,2.200000)
\uput[l](-2.117000,0){$|z|=1$}
\uput{10pt}[u](.995800,1.868000){$e^{i\phi}$}
\uput{10pt}[d](.995800,-1.868000){$e^{-i\phi}$}
\SpecialCoor\pscircle(0,0){2.117000}
\psline{->}(2.117000;0.000000)(3.146000;0.000000)
\psline{->}(2.117000;18.000000)(2.992000;18.000000)
\psline{->}(2.117000;36.000000)(2.646000;36.000000)
\psline{->}(2.117000;54.000000)(2.270000;54.000000)
\psline{->}(2.117000;72.000000)(1.938000;72.000000)
\psline{->}(2.117000;90.000000)(1.665000;90.000000)
\psline{->}(2.117000;108.000000)(1.452000;108.000000)
\psline{->}(2.117000;126.000000)(1.291000;126.000000)
\psline{->}(2.117000;144.000000)(1.180000;144.000000)
\psline{->}(2.117000;162.000000)(1.115000;162.000000)
\psline{->}(2.117000;180.000000)(1.093000;180.000000)
\psline{->}(2.117000;198.000000)(1.115000;198.000000)
\psline{->}(2.117000;216.000000)(1.180000;216.000000)
\psline{->}(2.117000;234.000000)(1.291000;234.000000)
\psline{->}(2.117000;252.000000)(1.452000;252.000000)
\psline{->}(2.117000;270.000000)(1.665000;270.000000)
\psline{->}(2.117000;288.000000)(1.938000;288.000000)
\psline{->}(2.117000;306.000000)(2.270000;306.000000)
\psline{->}(2.117000;324.000000)(2.646000;324.000000)
\psline{->}(2.117000;342.000000)(2.992000;342.000000)
\end{pspicture}
\caption{Gradient of $\Re S(z)$ on the unit circle}
\label{f6}
\end{center} \end{figure}
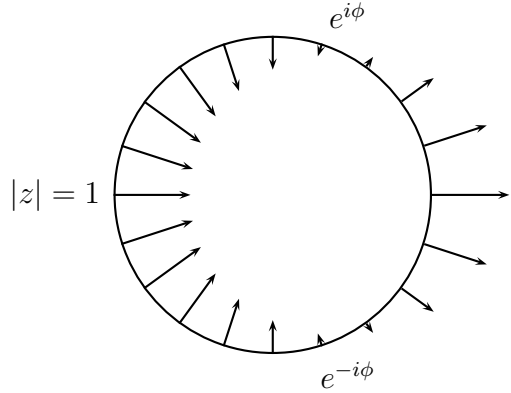

\begin{figure}[!h]
\centering
\scalebox{.8}{\includegraphics{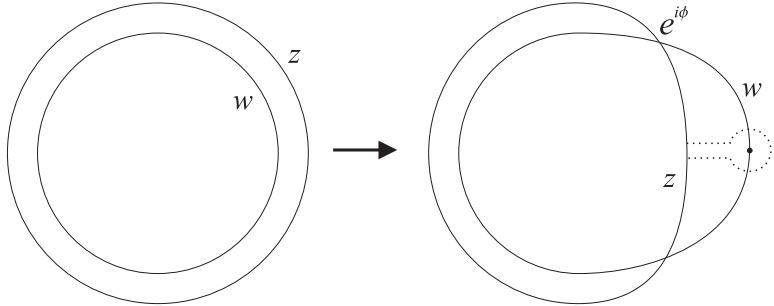}}
\caption{Deformation of the contours for the integral \eqref{Ib}}
\label{fsin}
\end{figure}

We now illustrate how this works on a particular example.
Namely, let us prove that as
$$
\frac{x}{\sqrt\xi},\frac{y}{\sqrt\xi} \to u\in (-2,2)\,,
$$
in such a way that $x-y$ remains fixed, we have
\begin{equation}
  \label{BtS}
\Ke{Bessel}(x,y;\xi) \to 
\Ke{sin}(x-y;\phi) \,,  
\quad \phi=\arccos(u/2) \,.
\end{equation}
We have
\begin{multline}
  \Ke{Bessel}(x,y;\xi) = \\ \frac1{(2\pi i)^2} \iint_{|z|>|w|}
  \frac{\exp\left(\sqrt{\xi}\left[S\left(z;\frac{x}{\sqrt\xi}\right)-
        S\left(w;\frac{y}{\sqrt\xi}\right) \right]\right)}{z-w} \,
  \frac{dz\, dw}{\sqrt{zw}}\,,  \label{Ib}
\end{multline}
where
$$
S(z;u)=z-z^{-1}-u \ln z \,.
$$
The critical points of the function $S(z,u)$ are the roots of
\begin{equation}
z S'(z;u)=z + z^{-1} - u 
\label{Sd}
\end{equation}
which for $u\in (-2,2)$ are clearly given by
$$
z=e^{\pm i \phi} \,.
$$
It is easy to see that the function $\Re S(z;u)$
vanishes on the unit circle $|z|=1$ and the direction 
of its gradient, plotted in Figure \ref{f6}, is easy
to infer from \eqref{Sd}. 
Thus, by the method of steepest descent we have to deform 
the contours of integration as shown in Figure \ref{fsin},
after which the integral, as explained above,
 will go to $0$ as $\xi\to\infty$.
However, as we deform the contours, we will pick up the residue
at $z=w$, as is also shown in Figure \ref{fsin}. It is obvious 
that the residue of the integrand in \eqref{Ib} at
$z=w$ is equal to $\dfrac{1}{2\pi i} \dfrac{dw}{w^{x-y+1}}$.
Thus, we conclude that
$$
\Ke{Bessel}(x,y;\xi) \to  \frac1{2\pi i}\int_{e^{-i\phi}}^{e^{i\phi}}
\frac{dw}{w^{x-y+1}} = \Ke{sin}(x-y;\phi)\,,
$$
as was to be shown.

\subsection{}

We now turn to the asymptotics on the edge of the limit
shape which for the Plancherel measure is described in 
Theorem \ref{TA}. Recall that near the edge $u=2$ of the limit
shape the distance between two consecutive ``downs'' is 
of order $N^{1/6}$ and hence also the probability to 
observe a ``down'' in any particular position is of
order $N^{-1/6}$, where $N$ is the size of the partition 
$\lambda$. For the Poissonized Plancherel measure $\Pb{PP}$
the expected size of $\lambda$ is $\xi$. From Theorem \ref{TA} 
we, therefore, expect that
\begin{equation}
\xi^{1/6} \, \Ke{Bessel}\left(2\sqrt \xi+x \xi^{1/6}, 
2\sqrt \xi+y \xi^{1/6}\right) \to \Ke{Airy}(x,y)\,,
\quad \xi\to\infty\,. \label{BtA}
\end{equation}
We will now use the steepest descent method 
to deduce \eqref{BtA} from \eqref{Ib}. The 
limit \eqref{BtA} is the key ingredient of 
the proofs of Theorem \ref{TA} given in \cite{BOO,J1}.
The difference between \eqref{BtA} and Theorem \ref{TA}
is a technical, but rather standard, collection of 
estimates that allows one to deduce the convergence
of joint distributions for Plancherel measures from
the convergence of correlation functions of their
Poissonizations.

\subsection{}

We now start the asymptotic analysis of the integral 
\eqref{Ib} in the $u=2$ regime. What happens in this
case is that the two nondegenerate critical points $e^{\pm i\phi}$
of $S(z,u)$ coalesce at $z=1$ to form one degenerate
critical point
\begin{equation}
S(z,2)=\tfrac13(z-1)^3+O((z-1)^4)\,.
\label{S1}
\end{equation}
Since $S(z,u)$ appears in the exponent in \eqref{Ib},
this starts looking like the Airy integral \eqref{IA},
and, indeed, this is precisely the source of the appearance of
the Airy function in \eqref{BtA}. In fact, this is 
the fundamental reason why the Airy function is as ubiquitous
in asymptotics as the sine function. The sine function 
appearance in the steepest descent asymptotics is usually
produced by two complex conjugate critical points. When the
parameters are tuned so that these two points 
coalesce, the Airy function enters the scene. This is all 
very standard and is at length discussed, for example, in 
\cite{W} and many other places.  

\begin{figure}[!h]
\centering
\scalebox{.7}{\includegraphics{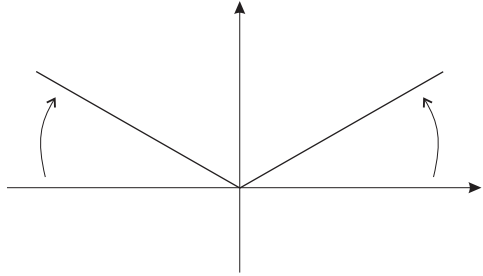}}
\caption{A better contour for the Airy integral \eqref{IA}}
\label{fcA}
\end{figure}

\begin{figure}[!h]
\centering
\scalebox{.45}[.3]{\includegraphics{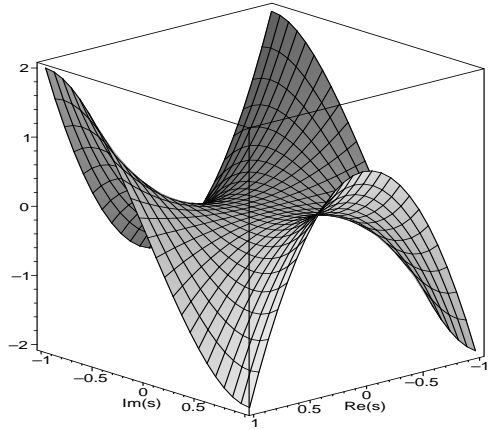}}
\caption{Plot of $\Re(i s^3)$}
\label{fsed}
\end{figure}

The integral \eqref{IA} converges, but only conditionally,
due to ever faster oscillation of the cosine function for
large $s$. A much better integral can be obtained from 
\eqref{IA} by shifting the contour of integration as 
shown in Figure \ref{fcA}. 
The contour in Figure \ref{fcA} goes along 2 of the 3 downslopes
 of the  
graph of $\Re(i s^3)$, plotted in Figure \ref{fsed} and 
commonly known as the monkey saddle. 
The function 
$e^{is^3/3}$ is very rapidly decaying on this contour. 

\subsection{}

\begin{figure}[!h]
\centering
\scalebox{.8}{\includegraphics{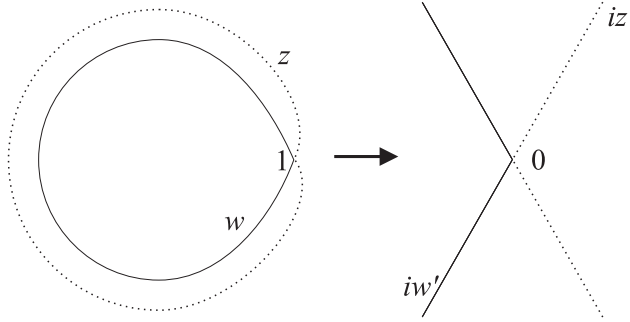}}
\caption{Contours for the edge ($u=2$) asymptotics}
\label{fai}
\end{figure}

It is clear that for $u=2$ all the action in the integral 
\eqref{Ib} is happening in the neighborhood of the 
critical point $z=w=1$, so we want to zoom in on this 
point. We also want the integration contour in $z$ to pass
through the downslopes of the graph of the real part of 
\eqref{S1}, while the $w$-contour should go along the 
ridges of the same graph  (because the of the opposite sign
of the $S\left(w;\frac{y}{\sqrt\xi}\right)$ term in \eqref{Ib}). 
We thus deform the 
integration as shown in the left half of Figure \ref{fai}
and introduce the following
new integration variables
\begin{equation}
z=1 + i \frac{z'}{\xi^{1/6}} \,, \quad 
w=1 + i \frac{w'}{\xi^{1/6}} \,.\label{zp}
\end{equation}
For large $\xi$, the integration contour in $z'$ and
$w'$ will look like the contour in the right half of
Figure \ref{fai}.
The scaling in \eqref{zp} is chosen so that to make the
exponent in \eqref{Ib}
$$
\sqrt \xi S(z,2)= - \frac{i}{3}\, (z')^3 + o(1) 
$$
stay finite as $\xi\to \infty$. 

More generally, the leading order expansion of $S(z,u)$
near $(z,u)=(1,2)$ with respect to both arguments is
given by 
$$
\sqrt \xi 
S\left(z,2+\frac{x}{\xi^{1/3}}\right) = 
- i x z' - \frac{i}{3}\, (z')^3  + o(1) \,, \quad \xi\to\infty \,.
$$
Plugging this into the integral \eqref{Ib} we obtain
the following
\begin{multline}
\xi^{1/6} \, \Ke{Bessel}\left(2\sqrt \xi+x \xi^{1/6}, 
2\sqrt \xi+y \xi^{1/6}\right) \to \\
\frac{1}{(2\pi)^2}
\iint \frac{e^{-i(z')^3/3-i z' x + i (w')^3/3 + i w' y}}{i(z'-w')}
\, dz' dw'\,,  \label{BtI}
\end{multline}
where the integration contour (rotated by $90^\circ$) is shown 
in the right half of Figure \ref{fai}.

\subsection{}

It remains to show that
\begin{multline}
  \label{IKA}
 \frac{1}{(2\pi)^2}
\iint \frac{e^{-i(z')^3/3-i z' x + i (w')^3/3 + i w' y}}{i(z'-w')}
\, dz' dw' = \\ \frac{\Ai(x) \, \Ai'(y) - \Ai'(x) \, \Ai(y)}{x-y} \,. 
\end{multline}
This is entirely parallel to the argument in Section \ref{sKB}
(see the more general discussion in Section 2.2.4  of \cite{O2}). 
We have
$$
e^{-i z' x + i w' y} =
\frac{i}{x-y} \, \left(
\frac{\partial}{\partial z'} + \frac{\partial}{\partial w'}
\right)  \, e^{-i z' x + i w' y} \,.
$$
Inserting this into the LHS of \eqref{IKA} and integrating
by parts we obtain 
\begin{multline*}
  \frac{1}{(2\pi)^2} \iint \frac{e^{-i(z')^3/3-i z' x + i (w')^3/3 + i
      w' y}}{i(z'-w')} \, dz' dw' =\\
 \frac{i}{(2\pi)^2\, (x-y)} \iint
  (z'+w') \, e^{-i\frac{(z')^3}3-i z' x + i \frac{(w')^3}3 + i w' y} \, dz'
  dw' \,.
\end{multline*}
The formula \eqref{IA} translates into 
$$
\frac{1}{2\pi} \int
  \, e^{- i (z')^3/3 - i z' x} \, 
  dz' = \Ai(x) \,,
$$
because our contour of integration in $z'$ (see Figure \ref{fai})
is the negative of the contour in Figure \ref{fcA}. It is also 
clear that 
$$
\frac{i}{2\pi} \int
  w' \, e^{i (w')^3/3 + i w' y} \, 
  dw' = \Ai'(y) \,.
$$
Thus, equation \eqref{IKA} follows.

\subsection{}

After finishing these computations, it may be useful to 
look at them from a more abstract point of view. 
Above, we tried to be absolutely concrete and explicit, but
in fact we used only very general properties of the 
correlation kernel \eqref{Ib}. Namely, all that was needed
was the information about the location of the 
critical points of the function $S$ in \eqref{Ib} and 
the residue of the integrand in \eqref{Ib} at $z=w$. 
The fact that the two critical were complex conjugate,
which was important for obtaining the sine function, is 
an automatic consequence of the fact the integrand is 
real for real values of the argument. 
Even less was required for the Airy function asymptotics,
the only essential ingredient in which was the coalescence
of the two critical points.

This means the same method is applicable to much a wider
variety of problems producing the same or equivalent
results. This phenomenon, when the final asymptotic result proves
to be insensitive to the fine details of the original problem,
is known as \emph{universality}. The simplest example of this
phenomenon may be the central limit theorem which says
that, under very general assumptions, sums of independent,
identically distributed, but otherwise arbitrary random 
variables converge, scaled and centered, to the normal
distribution. It is widely believed that
the universality principle applies very generally. From 
a mathematical perspective it is often quite mysterious why
this should be the case, but sometimes, like for the
central limit theorem, there is a transparent explanation.
Similarly, our exact formulas for the correlation functions
of the Schur process, combined with the steepest descent
argument, offer a simple explanation for the appearance of
the sine and Airy kernel.

\end{document}